\documentclass[12pt]{amsart} \usepackage{amssymb}
\usepackage{graphicx}

\theoremstyle{plain}
\newtheorem{theorem}{Theorem}
\newtheorem{lemma}[theorem]{Lemma}

\theoremstyle{definition}

\newtheorem{remark}[theorem]{Remark}

\newcommand{\ovl}{\overline}

\begin{document}

\author{G\'abor Kun}
\email{kungabor@cs.elte.hu}
\author{Jaroslav Ne\v{s}et\v{r}il}
\email{nesetril@kam.mff.cuni.cz}
\address{Department of Algebra and Number Theory, E\"otv\"os Lor\'and
University, 1117 P\'azm\'any P\'eter s\'et\'any 1/c}
\address{Department of Applied Mathematics (KAM) \\
   and Institute of Theoretical Computer Science (ITI)\\
   Charles University \\
   Malostransk\'{e} n\'{a}m 22 \\
   Praha\\}

\title{Forbidden Lifts\\(NP and CSP for combinatorists)}

\date{\today}
\thanks{Part of this work was supported by ITI and DIMATIA of Charles
University Prague under grant 1M0021620808, by OTKA Grant no. T043671
and also by Isaac Newton Institute (INI) Cambridge.}

\maketitle
\begin{abstract}
We present a definition of the class NP in combinatorial  context as
the set of languages of structures defined by finitely many
forbidden lifted substructures. We apply this to special
syntactically defined subclasses and show how they correspond to
naturally defined (and intensively studied) combinatorial problems.
We show that some types of combinatorial problems like edge
colorings and graph decompositions express the full computational
power of the class NP. We then characterize Constraint Satisfaction
Problems (i.e. $H$-coloring problems) which are expressible by
finitely many forbidden lifted substructures. This greatly
simplifies and generalizes the earlier attempts to characterize this
problem. As a corollary of this approach we perhaps find a proper
setting of Feder and Vardi analysis of CSP languages within the
class MMSNP.
\end{abstract}

\section{Introduction}
Think of a $3$-colorability of a graph $G = (V,E)$. This is a well known
hard problem and there is a multiple evidence for this: concrete
instances of the problem are difficult to solve (if you want a
non-trivial example consider Kneser graphs; \cite{UBM}), there is
an abundance of minimal graphs which are not $3$-colorable (these
are called $4$-critical graphs, see e.g. \cite{JT}) and in the
full generality (and even for important ``small" subclasses such
as $4$-regular graphs or planar graphs) the problem is a canonical
NP-complete problem.

Yet the problem has an easy formulation. A $3$-coloring is simple
to formulate even at the kindergarten level. This is in a sharp
contrast with the usual definition of the class NP by means of
polynomially bounded non-deterministic computations. Fagin
\cite{FA} gave a concise description of the class NP by means of
logic: NP languages are just languages accepted by an Existential
Second Order (ESO) formula of the form

\[
\exists P \Psi(S,P),
\]
where $S$ is the set of input relations, $P$ is a set of existential
relations, the proof for the membership in the class, and $\Psi$ is
a first-order formula without existential quantifiers. This
definition of NP inspired a sequence of related investigations (see
e.g. \cite{DKL,IMM,V} and these {\it descriptive complexity} results
established that most major complexity classes can be characterized
in terms of logical definability of finite structures. Particularly
this led Feder and Vardi \cite{FV} to their seminal reduction of
{\it Constraint Satisfaction Problems} (shortly CSP) to the so
called MMSNP ({\it Monotone Monadic Strict Nondeterministic
Polynomial}) problems which also nicely link MMSNP to the class NP
in computational sense. This will be explained in some detail in
Section 3 which presents one of the main motivations of this paper.
Inspired by these results we would like to ask an even simpler
question:

\vspace{5mm}
Can one express the computational power of the class NP by combinatorial
means?
\vspace{5mm}

From the combinatorial point of view there is a standard way how to
approach (and sometimes to solve) a monotone property $P$: one
investigates those structures without the property $P$ which are
{\it critical}, (or {\it minimal}) without $P$. One proceeds as
follows: denote by $\mathcal F$ the class of all  critical
structures and define the class ${\rm Forb}(\mathcal F)$ of all
structures which do not ``contain'' any $F \in \mathcal F$. The
class ${\rm Forb}(\mathcal F)$ is the class of all structures not
containing any of the critical substructures and thus it is easy to
see that ${\rm Forb}(\mathcal F)$ coincides with the class of
structures with the property $P$. Of course in most cases the class
$\mathcal F$ is infinite yet a structural result about it may shed
some light on property $P$. For example this is the case with
$3$-colorability of graphs where $4$-critical graphs were (and are)
studied thoroughly (historically mostly in relationship to the Four
Color Conjecture).

Of particular interest (and as the extremal case in our setting)
are  those monotone properties $P$ of structures which can be
described by finitely many forbidden substructures. It has been
proved in a sequence of papers \cite{atserias,ros} that a
homomorphism monotone problem is {\it First Order} (shortly FO)
definable if and only if it is {\it positively} FO definable
(shortly FO+ definable), i.e. the formula does not contain any
negations (and so implications and inequalitites), and thus
alternatively defined as  ${\rm Forb}(\mathcal F)$ for a finite
set $\mathcal F$ of structures. Although FO-definability is not a
rare fact (and extremely useful in database theory), still
FO-definability cannot express most combinatorial problems
(compare \cite{NT},\cite{atserias} which  characterize all CSP
which are FO-definable; see also Theorem \ref{FNT}). Thus it may
seem to be surprising that the classes of relational structures
defined by ESO formulas (i.e. the whole class NP) corresponds
exactly to those canonical {\it lifts} of structures which are
defined by a finite set of forbidden substructures. Shortly,
finitely many forbidden lifts determine any language in NP. This
is being made precise in Section 3. Here, let us just briefly
illustrate this by our example of $3$-colorability. Instead of a
graph $G = (V, E)$ we consider the graph $G$ together with three
unary relations $C_1, C_2, C_3$ which $\it cover$ the vertex set
$V$; this structure will be denoted by $G'$ and called a {\it
lift} of $G$ ($G'$ has one binary and three unary relations).
There are $3$ {\it forbidden substructures} or {\it patterns}: For
each $i = 1,2,3$ the graph $K_2$ together with cover $C_i =
\{1,2\}$ and $C_j = \emptyset$ for $j \neq i$ form pattern ${\bf
F'_i}$ (where the signature of ${\bf F'_i}$ contains one binary and
three unary relations). The language of all $3$-colorable graphs then
corresponds just to the language $\Phi({\rm Forb}({\bf F'_1},{\bf
F'_2},{\bf F'_3}))$ where $\Phi$ is the forgetful functor which
transforms $G'$ to $G$ and the language of $3$-colorable graphs is
just the language of the class satisfying formula $\exists
{G'}({G'} \in {\rm Forb}({\bf F'_1},{\bf F'_2},{\bf F'_3}))$. This
{\it extended language} (of structures $G'$) of course just
expresses the membership of  $3$-colorability to the class NP.
There is more than this that meets the eye. This scheme fits
nicely into the mainstream combinatorial and combinatorial
complexity research. Building upon Feder-Vardi classification of
MMSNP we isolate (in Theorems \ref{NP}, \ref{inj}, \ref{full})
three computationally equivalent formulations of NP class:

\begin{enumerate}

\item
By means of shadows of forbidden homomorphisms of relational lifts (the corresponding category is denoted by $Rel^{cov}(\Delta, \Delta')$),

\item
By means of shadows of forbidden injections (monomorphisms) of
monadic lifts (the corresponding category will be denoted by
$Rel^{cov}_{inj}(\Delta, \Delta')$),

\item
By means of shadows of forbidden full homomorphisms of monadic
lifts (the corresponding category will be denoted by
$Rel^{cov}_{full}(\Delta, \Delta')$).
\end{enumerate}

Our results imply that each of these approaches includes the whole
class NP. It is interesting to note how nicely these categories
fit to the combinatorial common sense about the difficulty of
problems: On the one side the problems in CSP correspond and
generalize ordinary (vertex) coloring problems. One expects a {\it
dichotomy} here: every CSP problem should be either polynomial or
NP-complete (as conjectured in \cite{FV} and probabilistically
verified in \cite{LN}). On the other side the above formulations (1),
(2), (3) model the whole class NP and thus we cannot expect dichotomy
there (by a celebrated result of Ladner \cite{Lad}). But this is
in accordance with the combinatorial meaning of these classes: the
formulation (1) expresses coloring of edges, triples etc. and thus it
involves problems in Ramsey theory \cite{GRS,N}. The formulation (2)
may express vertex coloring of classes with restricted degrees of
vertices \cite{K,HN}. The formulation (3) relates to vertex colorings
with a given pattern among classes which appears in many graph
decomposition techniques (for example in the solution of the
Perfect Graph Conjecture \cite{CST}). The point of view of
forbidden partitions (in the language of graphs and matrices) is
taken for example in \cite{HELL}. This clear difference between
combinatorial interpretations of syntactic restrictions on
formulas expressing the computational power of NP is one of the
pleasant consequences of our approach.

At this point we should add one more remark. We of course do not
only claim that every problem in NP can be polynomially reduced to
a problem in one of these classes. This would only mean that each
of these classes contains an NP-complete problem. What we claim is
that these classes have the {\it computational power} of the whole
of NP, i.e. these classes are {\it computationally equivalent} to
all problems in NP. More precisely, to each language $L$ in NP
there exists a language $M$ in any of these three classes such
that $M$ is {\it polynomially equivalent} to $L$, i.e. there exist
{\it polynomial reductions} of $L$ to $M$ and $M$ to $L$.


Having finitely many forbidden patterns (i.e. forbidden
substructures) for a class of structures $\mathcal K$ we are
naturally led to the question whether $\mathcal K$ is the class
determined by a finite set of templates, or in other words by the
existence of homomorphisms to particular structures. In technical
terms (see e.g. \cite{HN, FV}) this amounts to the question
whether $\mathcal K$ is an instance of a {\it Constraint
Satisfaction Problem} (shortly CSP). On the other hand finitely
many forbidden patterns lead to the question whether the class
$\mathcal K$ is not defined by a {\it finite duality}. This scheme
for combinatorial problems goes back to \cite{NP}, see e.g.
\cite{HN} and it was studied in situations as diverse as {\it
bounded tree width dualities} \cite{HNZ}, {\it duality of linear
programming} \cite{HN} and {\it classes with bounded expansion}
\cite{NPOM}. Here we completely characterize (using results of
\cite{NT}) shadows of finitary dualities in the case where the
extension of the language is monadic, i.e. it consists of unary
relations (as is the above case of $3$-coloring), see Theorem
\ref{main}. These general results can be used in the investigation
of the class MMSNP (to be defined in Section 3). Feder and Vardi
introduced this class as a fragment of SNP in \cite{FV}. They
proved that the class MMSNP is randomly polynomially equivalent to
the class of finite union of CSP languages. This was later
derandomized by the first author proving that the classes MMSNP
and CSP are computationally equivalent \cite{CSPvsMMSNP}. We will
examine these classes from the viewpoint of descriptive complexity
theory: Any finite union of CSP languages belongs to MMSNP. But
the converse does not hold. Consider for example the language of
triangle free graphs: this is an MMSNP language which is not a
finite union of CSP languages. Madelaine and Stewart introduced
the class of {\it Forbidden Pattern Problems} (FP) as an
equivalent combinatorial version of MMSNP \cite{MStuart},
\cite{M}. They gave an effective, yet lengthy process to decide
whether an MMSNP language is a CSP language. We give a short and
easy procedure to decide whether an MMSNP language is a finite
union of CSP languages, and we show that these are exactly those
languages defined by forbidden patterns not containing any cycle.
This simplicity is possible by translation and generalization of
the Feder-Vardi proof of the computational equivalence of finite
union of CSP's and MMSNP in the context of category theoretical
language of duality.

The paper is organized as follows: In Section 2 we review the
basic notions and previous work related to finite structures.
Particularly we state two our basic tools: the characterization of
{\it finite dualities} \cite{NT, FNT} and a combinatorial
classique,  the {\it sparse incomparability lemma}. It is here
where we introduce two our basic notions of {\it lifts} and {\it
shadows}. The interplay of corresponding classes (categories) is a
central theme of this paper. In Section 3 we introduce the
relevant notions  of descriptive complexity (mostly taken from
\cite{FV}) and relate it to our approach. We prove that the class
NP is polynomially equivalent  with classes of structures
characterized by finitely many forbidden lifts (this is proved in
three different categories, see Theorems \ref{NP}, \ref{inj} and
\ref{full}). In Section 4 we study the relationship of lifts and
shadows abstractly from the point of view of dualities. Theorem
\ref{main} enables us to prove the characterization of shadows of
finite dualities  (called {\it lifted dualities}) in lifts and
shadows. This, as a corollary, proves the main result of
\cite{MStuart}. In Section 5 we return to Feder-Vardi setting and
indicate how the polynomial equivalence of classes MMSNP and
finite unions of CSP problems  emerges naturally in our
combinatorial-categorical context.

For more complicated (i.e. nonmonadic) lifts we (of course) have
partial results only. Perhaps the next case is that of covering
equivalences. This we are still able to handle with our methods
and we characterize all CSP languages in this class. But we postpone
this to another occasion.

\section{Categories of Finite Structures}

For a relational symbol $R$ and relational structure $\bf A$ let $A
= X(\bf A)$ denote the universe of $\bf A$ and let $R(\bf A)$ denote
the relation set of tuples of $\bf A$ which belong to $R$. Let
$\Delta$ denote the {\it signature} (type) of relational symbols,
and let $Rel(\Delta)$ denote the class of all relational structures
with signature $\Delta$. We will often work with two (fixed)
signatures, $\Delta$ and $\Delta \cup \Delta'$ (the signatures
$\Delta$ and $\Delta'$ are always supposed to be disjoint). For
convenience we denote structures in $Rel(\Delta)$ by $\bf A,\bf B$
etc. and structures in $Rel(\Delta \cup \Delta')$ by $\bf A',\bf B'$
etc. We shall denote $Rel(\Delta \cup \Delta')$ by $Rel(\Delta,
\Delta')$. The classes $Rel(\Delta)$ and $Rel(\Delta, \Delta')$ will
be considered as categories endowed with all homomorphisms. Recall,
that a homomorphism is a mapping which preserves  all relations.
Just to be explicit, for relational structures ${\bf A,\bf B} \in
Rel(\Delta)$ a mapping $f: X({\bf A}) \longrightarrow X(\bf B)$ is a
{\it homomorphism} ${\bf A} \longrightarrow \bf B$ if for every
relational symbol $R \in \Delta$ and for every tuple $(x_1, \ldots,
x_t) \in R(\bf A)$ we have $(f(x_1), \ldots, f(x_t)) \in R(\bf B)$.
More generally we will use the notation ${\bf A} \longrightarrow \bf
B$ for morphisms when working in other categories. (These will be
categories of relational structures, where the morphisms will be
either the injective or the full homomorphism, respectively.)
Similarly we define homomorphisms for the class $Rel(\Delta,
\Delta')$. The interplay of categories $Rel(\Delta, \Delta')$ and
$Rel(\Delta)$ is one of the central themes of this paper. Towards
this end we define the following notions: Let $\Phi: Rel(\Delta,
\Delta') \rightarrow Rel(\Delta)$ denote the natural {\it forgetful
functor} that ``forgets'' relations in $\Delta'$. Explicitly, for a
structure ${\bf A'} \in Rel(\Delta, \Delta')$ we denote by $\Phi(\bf
A')$ the corresponding structure ${\bf A} \in Rel(\Delta)$ defined
by $X({\bf A'}) = X(\bf A)$, $R({\bf A'}) = R(\bf A)$ for every $R
\in \Delta$. For a homomorphism $f: {\bf A}' \longrightarrow \bf B'$
we put $\Phi(f) = f$. The mapping $f$ is of course also a
homomorphism $\Phi({\bf A'}) \longrightarrow \Phi(\bf B')$. This is
expressed by the following diagram.

\begin{figure}[!h]
\begin{center}
\includegraphics[width=4.5cm]{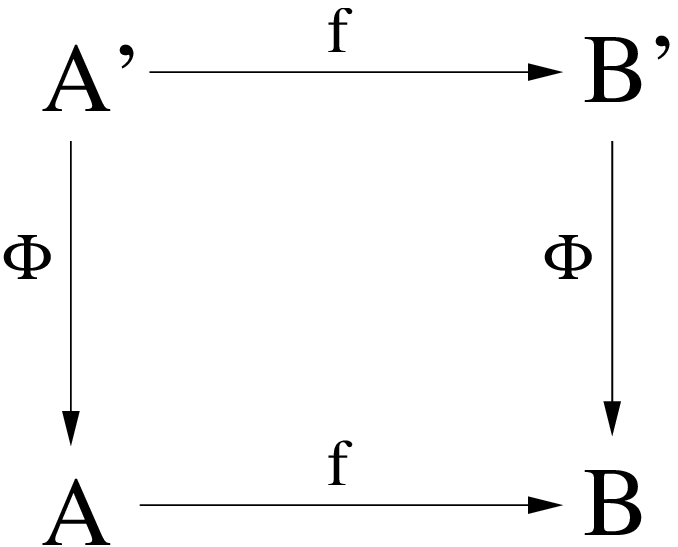}
\end{center}
\end{figure}

\vspace{10mm}

These object-transformations call for a special terminology: For
${\bf A}' \in Rel(\Delta, \Delta')$ we call $\Phi(\bf A') = \bf A$
the {\it shadow} of $\bf A'$. Any $\bf A'$ with $\Phi(\bf A') = \bf
A$ is called a {\it lift} of $\bf A$. The analogous terminology is
used for subclasses $\mathcal C$ of $Rel(\Delta, \Delta')$ and
$Rel(\Delta)$. (Thus, for example, for a subclass $\mathcal{C}
\subseteq Rel(\Delta, \Delta'),  \Phi(\mathcal C)$ is the class of
all shadows of all structures in the  class $\mathcal C$.) The
following special subclass of $Rel(\Delta, \Delta')$ will be
important: denote by $Rel^{cov}(\Delta, \Delta')$ the class of all
structures in $Rel(\Delta, \Delta')$ where we assume that all
relations in $\Delta'$ have the same arity, say $r$, and that all
the $r$-tuples of an object are contained by some relation in
$\Delta'$. The category $Rel^{cov}(\Delta, \Delta')$ is briefly
called {\it covering} or $r$-{\it covering category}. Note that the
class ${Rel^{cov}(\Delta, \Delta')}$ is closed under surjective
homomorphisms. We will work with two other similar pairs of
categories. We denote by $Rel_{inj}(\Delta)$ and
$Rel_{full}(\Delta)$ the categories where the objects are again the
relational structures of type $\Delta$, but the morphisms are the
injective and full homomorphisms, respectively. We call a mapping a
{\it full homomorphism} if it is relation and non-relation
preserving, too. Such mappings have very easy structure, as every
full homomorphism which is onto is a retraction. We denote by
$Rel^{cov}_{inj}(\Delta,\Delta')$ and
$Rel^{cov}_{full}(\Delta,\Delta')$ the subclasses containing the
same class of objects as $Rel^{cov}(\Delta,\Delta')$. We only will
use these notions in the case when $\Delta'$ contains monadic
relations.

Let $\mathcal F'$ be a finite set of structures in the category
$\mathcal C$ (one of the above categories). By ${\rm
Forb}(\mathcal F')$ we denote the class of all structures ${\bf
A}' \in {\mathcal C}$ satisfying ${\bf F}' \not\longrightarrow \bf
A'$ for every ${\bf F}' \in \mathcal F'$. (This class is sometimes
and perhaps more efficiently denoted by ${\mathcal F'}
\not\rightarrow$.) Similarly (well, dually), for the finite set of
structures $\mathcal D'$ in $\mathcal C$ we denote by
$CSP(\mathcal D')$ the class of all structures ${\bf A}' \in
\mathcal C$ satisfying ${\bf A}' \longrightarrow {\bf D}'$ for
some ${\bf D}' \in \mathcal D'$. (This is sometimes denoted  by
$\rightarrow \mathcal D$.) Now suppose that the classes ${\rm
Forb}(\mathcal F')$ and $CSP(\mathcal D')$ are equal. Then we say
that the pair $(\mathcal {F', D'})$ is a {\it finite duality} in
$\mathcal C$. Explicitly, a finite duality  means that the
following equivalence holds for every structure ${\bf A}' \in
\mathcal C$:

\noindent ${\bf F}' \not\longrightarrow {\bf A}' \text{ for every }
{\bf F}' \in \mathcal F' \text{ iff } {\bf A}'  \longrightarrow {\bf
D}' \text{ for some } {\bf D}' \in \mathcal D'.$


One more definition is needed. In dualities (as well as in most of
this paper) we are interested in the existence of a homomorphism
(every CSP can be expressed by the existence of a homomorphism to
a template; see \cite{FV},\cite{HN}) . Consequently we can also use
the language of partially ordered sets and consider the {\it
homomorphism order} $\mathcal C_{\Delta}$ defined on the class of
all structures with signature $\Delta$: we define the order $\leq$
by putting ${\bf A} \leq {\bf B}$ iff there is a homomorphism
${\bf A} \longrightarrow {\bf B}$. The ordering $\leq$ is clearly
a quasiorder but this becomes a partial order if we either
factorize $\mathcal C_{\Delta}$ by the homomorphism equivalence
or, perhaps preferably, if we restrict $\mathcal C_{\Delta}$ to
non-isomorphic $\it core$ structures. We say that $\bf A$ is  $\it
core$ if every homomorphism ${\bf A} \longrightarrow {\bf A}$ is
an automorphism. Every finite structure $\bf A$ contains (up to an
isomorphism) a uniquely determined core substructure homomorphically
equivalent to $\bf A$, see \cite{NT, HN}. The following result was
recently proved in \cite{FNT} as a generalization of \cite{NT}. It
characterizes finite dualities of finite structures., i.e. in the
category $Rel(\Delta)$.

\begin{theorem}
\label{FNT}
For every signature $\Delta$ and for every finite set $\mathcal F$ of
(relational) forests there exists (up to a homomorphism equivalence) a
uniquely determined set $\mathcal D$ of structures
such that $(\mathcal {F, D})$ forms a finite duality.
Up to a homomorphism equivalence there are no other finite dualities.
\end{theorem}

We did not define what is a forest (see \cite{NT,FNT}). For the sake
of completeness let us say that a {\it forest} is a structure not
containing any cycle. And a {\it cycle} in a structure $\bf A$ is
either a sequence of distinct points and distinct tuples $x_0, r_1,
x_1, \ldots, r_t, x_t = x_0$ where each tuple $r_i$ belongs to one
of the relations $R(\bf A)$ and each $x_i$ is a coordinate of $r_i$
and $r_{i+1}$, or, in the degenerated case $t = 1$ a relational
tuple with at least one multiple coordinate. The {\it length} of the
cycle is the integer $t$ in the first case and $1$ in the second
case. Finally the {\it girth} of a structure $\bf A$ is the shortest
length of a cycle in $\bf A$ (if it exists; otherwise it is a
forest).

The study of homomorphism properties of structures not containing
short cycles (i.e. with large girth) is a combinatorial problem
studied intensively. The following result has  proved particularly
useful in various applications. It is often called the {\it Sparse
Incomparability Lemma}:

\begin{lemma}
\label{sparse}

Let $k, \ell$ be positive integers and let $\bf A$ be a structure. Then there exists a structure $\bf B$ with the following properties:

\begin{enumerate}
\item There exists a homomorphism $f: {\bf B} \longrightarrow {\bf A}$;

\item For every structure $\bf C$ with at most $k$ points the following holds:
there exists a homomorphism $ {\bf A} \longrightarrow {\bf C}$ if and only if
there exists a homomorphism $ {\bf B} \longrightarrow {\bf C}$;

\item $\bf B$ has girth $\geq \ell$.
\end{enumerate}
\end{lemma}

\begin{figure}[!h]
\begin{center}
\includegraphics[width=3.0cm]{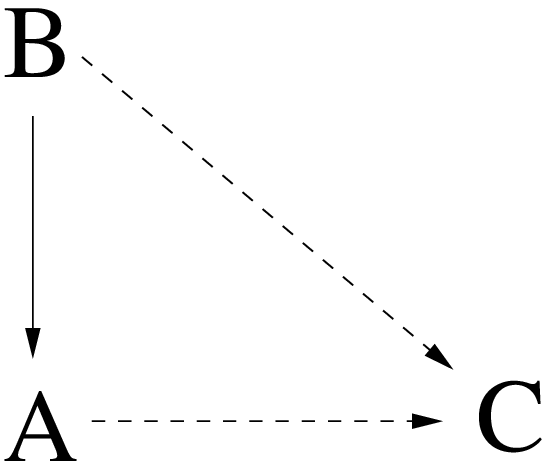}
\end{center}
\end{figure}

\vspace{6mm}

This result was proved in \cite{NR, NZ} (see also \cite{HN}) by
probabilistic methods. In fact in \cite{NR, NZ} it was proved for
graphs only but the proof is the same for finite relational
structures. Of particular interest in this context is the question
whether there exists an explicit construction of the structure $\bf
B$. This is indeed possible: in the case of binary relations
(digraphs) this was done in \cite{MN} and for general relational
structures in \cite{CSPvsMMSNP}.

\section{NP by means of finitely many forbidden lifts}

There is  a standard connection between formulae  and existence of
homomorphisms.
This goes back to \cite{CHAM} and it can be formulated as follows:

To every structure $\bf A$ in $Rel(\Delta)$ we
associate the {\it canonical conjunctive existential formula}
$\varphi_{\bf A}$ as the conjunction of the atoms $R_{\bf
A}(\ovl{x})$, where $R \in \Delta$ preceded by existential
quantification of all elements of $\bf A$. Clearly this process
may be reversed and thus there is a one-to-one correspondence
between canonical conjunctive existential formulae and structures.
It is then obvious that the following holds:

There is a homomorphism ${\bf A} \longrightarrow {\bf B}$ if and only if
${\bf B} \models \varphi_{\bf A}$.

Following Fagin \cite{FA}, the class SNP consists of all problems
expressible by an existential second-order formula with a universal
first-order part. The class of problems expressible by an
existential second-order formula is exactly the class NP when
restricted to languages of finite structures. So the class SNP is
computationally equivalent to NP. The input of any problem in SNP is
a relational structure $\bf A$ of signature $\Delta$ with base set
$A = X(\bf A)$ and $\Pi$ is a set of relations on the same base set
$A$. In this situation it is customary to call the second order
relations  $\Pi$ {\it proof}. Let us be more specific (see
\cite{FV}). Every language (problem) $L$ in SNP may be  equivalently
described by a formula of the form

\[
\exists \Pi \forall \ovl{x} \bigwedge_i \neg \big( \alpha_i \wedge
\beta_i \wedge \varepsilon_i \big),
\]
where
\begin{enumerate}

\item
$\alpha_i$ is a conjunction of atoms or negated atoms involving variables and input relations (i.e. of the form $R(\ovl x)$ and $\neg R(\ovl x)$ for a relational symbol $R$ and $\ovl x$ a tuple of elements of $X$),

\item
$\beta_i$ is a conjunction of atoms and negated atoms involving variables and existential (proof)
relations (i.e. of the form $P(\ovl x)$ and $\neg P(\ovl x)$ for $P \in \Pi$
and $\ovl x$ a tuple of elements of $X$) and

\item
$\varepsilon_i$ is the conjunction of atoms involving variables and
inequalities (i.e. of form $x \neq y$).
\end{enumerate}

A formula of this type is called a {\it canonical formula} of the language $L$ in SNP.
It will be denoted by $\varphi_{L}$.

{\bf Example}: Consider the following language of digraphs (i.e.
relational structures, where the signature contains one single
binary relation $E$) defined by the following ESO formula:

\noindent $\exists P_1 \exists P_2  \forall x_1, x_2, x_3,y
\bigwedge_k \neg \big[ (P_k(x_1) \wedge P_k(x_2) \wedge P_k(x_3))
\wedge (E(x_1, x_2) \wedge E(x_1, x_3) \wedge E(x_2, x_3)) \wedge
(x_1 \neq x_2 \wedge  x_1 \neq x_3 \wedge x_2 \neq x_3) \big] \wedge
\big[ \neg( \neg P_1(y) \wedge \neg P_2(y))\big]$.

This formula corresponds to the language of all binary relations
whose base set can be covered by two sets in such a way that none of
these sets contains linearly ordered set with $3$ elements. If we in
addition postulate that the relation $E$ is symmetric then these are
just graphs which can be vertex partitioned into two triangle free
graphs. Following \cite{FV} one can also define three important
syntactically restricted subclasses of SNP.

We say that a canonical formula is {\it monotone} if there are no negations in
the $\alpha_i$'s. This implies that more relations lead to less satisfiable formulae.
The canonical formula is {\it monadic} if the relations in $\Pi$ are all
monadic (which means that all proof relations are unary).
The canonical formula is said to be  {\it without inequality} if it can described by a
canonical formula which does not contain $\varepsilon_i$.

Feder and Vardi \cite{FV} have proved that the three subclasses of
SNP defined by formulae with any two of these syntactical
restrictions still have the full computational power of the class
NP.

\begin{theorem}\cite{FV}
\label{123}
\begin{enumerate}

\item Every problem in NP has a polynomially equivalent problem in
monotone SNP without inequality. Moreover, we may assume that the
existential relations are at most binary.

\item Every problem in NP has a polynomially equivalent problem in
monotone, monadic SNP.

\item Every problem in NP has a polynomially equivalent problem in
monadic SNP without inequality.

\end{enumerate}
\end{theorem}

(The claim that we may restrict to binary relations in (1) is not
stated in \cite{FV} but it is clear from the proof.)
The class
with all the three restrictions is denoted by MMSNP ({\it Monotone
Monadic Strict Nondeterministic Polynomial}). We deal with this
class in Section 6.

In this paper we will reformulate Theorem \ref{123} in our
combinatorial category lift/shadow setting. This will be done in
Theorem \ref{NP} for item $(1)$, in Theorem \ref{inj} for item $(2)$
and in Theorem \ref{full} for item $(3)$. First, we introduce the
following: we say that the formula is {\it primitive} if for every
clause $\big( \alpha_i \wedge \beta_i \wedge \varepsilon_i \big)$,
every variables $x_1, \dots ,x_r$ occurring in it and every
existential relation $P \in \Pi$ of arity $r$ either the atom
$P(x_1, \dots ,x_r)$ or its negation is an atom of the clause. We
need the following technical lemma.

\begin{lemma}~\label{tech}
Every language in SNP can be described by a primitive formula.
Moreover, if the original formula satisfies some of the
restrictions (i.e. if it is either monotone or monadic or without
inequality) then so does the primitive formula.
\end{lemma}
\begin{proof}
Consider the language $L$ and the canonical formula defining $L$:

\noindent $\exists P \forall \ovl{x} \bigwedge_i \neg \big( \alpha_i
\wedge \beta_i \wedge \varepsilon_i \big)$. We modify the formula so
that for every proof relation $R$ of arity $r$ and variables $x_1,
\dots ,x_r \in S$ appearing in $\alpha_i$ or $\beta_i$ either
$R(x_1, \dots ,x_r)$ or $\neg R(x_1, \dots ,x_r)$ is in the
appropriate conjunct. In order to have such a formula we can replace
$\neg \Big(\alpha_i \wedge \beta_i \wedge \varepsilon_i \Big)$ by
$\neg \big( \alpha_i \wedge \beta_i \wedge \varepsilon_i \wedge
R(x_1, \dots ,x_r) \big) \wedge \neg \big( \alpha_i \wedge \beta_i
\wedge \varepsilon_i \wedge \neg R(x_1, \dots ,x_r) \big)$, this is
equivalent to the original formula. The repetition of this process
will terminate in finitely many steps, and it gives an appropriate
formula.
\end{proof}

\vspace{5mm}

\begin{theorem}~\label{NP}
For every language $L \in NP$ there exist relational types
$\Delta,\Delta'$ and a finite set
$\mathcal{F'}$ of structures in $Rel^{cov}(\Delta, \Delta')$ such that $L$ is
computationally equivalent to $\Phi(Forb(\mathcal{F'}))$. Moreover, we may
assume that the relations in $\Delta'$ are at most binary.
\end{theorem}

This theorem presents an equivalent form of item $(1)$ of
Theorem~\ref{123} by means of homomorphisms and classes
$Forb(\mathcal{F'})$. It is interesting to express other conditions
$(2)$, $(3)$ of Theorem \ref{123} by means of homomorphisms and
classes $Forb(\mathcal{F'})$. These two other versions are stated
below as Theorems \ref{inj} and \ref{full}.

\begin{proof}

Consider a language $L$ and the canonical formula $\varphi_L$
(showing that it is monotone SNP without inequality). The
construction of $\mathcal{F'}$ consists of two steps. In the first
step we enforce technical conditions on the formula.

\noindent
{\bf Step 1.}

\noindent We need the technical assumption that all proof relations
in $\Pi$ have the same (at most binary) arity and the formula is
primitive. The first condition can be achieved by exchanging
relational symbols not of maximal arity by new relational symbols of
maximal arity (binary would suffice). We can proceed as follows. In
every clause of the formula we put new (free) different variables
into the new entries in $\beta_i$, and we increase the number of
variables in $\ovl{x}$, too. This new formula is equivalent to the
original one. An evaluation satisfies to the new formula exactly iff
its restriction to the original variables satisfies the original
formula. By Lemma~\ref{tech} we may also assume that the new formula
is primitive. In the following we denote this formula by
$\varphi_L$.

In the second step we define lifts.

\noindent
{\bf Step 2.}

\noindent
The type $\Delta'$ will contain $2^{|\Pi|}$ relational
symbols corresponding to the $2^{|\Pi|}$ possibilities for a
subset of proof relations indicating possibilities in which a
tuple can be.
The pattern ${\bf F}_i'$ will correspond to the clause $\alpha_i
\wedge \beta_i$. The base set of each structure ${\bf F}_i'$ is the
set of variables in the clause $\alpha_i \wedge \beta_i$. A tuple
$\ovl{t}$ of variables is in a relation $R$ (of type $\Delta$) if
the atom $R(\ovl{t})$ appears in $\alpha_i$. Every tuple $\ovl{t}$
in ${\bf F}_i'$ (of appropriate arity) is in exactly one relation
from $\Delta'$, this is the relation corresponding to the subset of
all existential relations $P \in \Pi$ such that the atom
$P(\ovl{t})$ appears in $\beta_i$. Let $\mathcal{F'}$ be the set of
all lifts ${\bf F}_i'$. These may be disconnected, although we may
work with their connected components, see Remark~\ref{Victor}.

We prove that for a structure ${\bf A} \in Rel(\Delta)$ the formula
$\varphi_L$ is satisfiable iff there is a lifted structure ${\bf A'}
\in Rel(\Delta, \Delta')$ such that no ${\bf F}_i' \in \mathcal{F'}$
maps to $\bf A'$.

Suppose that ${\bf A} \in \Phi(Forb(\mathcal{F'}))$, i.e. there is a
lift ${\bf A'} \in Forb(\mathcal{F'})$. We may suppose that every
tuple of ${\bf A'}$ is in exactly one $\Delta'$ relation. This
single $\Delta'$ relation corresponds to a subset of relational
symbols in $\Pi$. For every relational symbol $P \in \Pi$ define
$P({\bf A})$ to be the set of tuples $\ovl{t}$ such that the
relation in $\Delta'$ containing $\ovl{t}$ corresponds to a subset
containing $P$. Denote this structure of type $(\Delta, \Pi)$ by
${\bf A''}$. We show that these relations prove ${\bf A} \models
\varphi_L$. We have to prove that ${\bf A''} \models \neg (\alpha_i
\wedge \beta_i)$ holds for every clause. Consider the corresponding
forbidden lift ${\bf F'}_i$. We know that ${\bf F'}_i \rightarrow
{\bf A'}$, which yields ${\bf A''} \models \neg (\alpha_i \wedge
\beta_i)$.

Secondly suppose that ${\bf A} \models \varphi_L$. We can correspond
to the proof relations on ${\bf A}$ a $\Delta'$ covering of the
(binary) tuples in ${\bf A}$, where every tuple is covered exactly
once. This lift ${\bf A'}$ shows that ${\bf A} \in
\Phi(Forb(\mathcal{F'}))$. Consider a forbidden lift ${\bf F'}_i \in
\mathcal{F'}$. We know that the proof relations of type $\Pi$
satisfy the formula $\neg (\alpha_i \wedge \beta_i)$, hence ${\bf
F'}_i \nrightarrow {\bf A}'$.
\end{proof}

\begin{remark}~\label{Victor}
Consider the languages $\Phi(Forb(\mathcal{F'}))$ and
$\Phi(Forb(\mathcal{G'}))$. Their union is exactly the language
$\Phi(Forb(\mathcal{H'}))$, where $\mathcal{H'}=\{ {\bf F'}\cup
{\bf G'}: {\bf F'} \in Forb(\mathcal{F'}), {\bf G'} \in
Forb(\mathcal{G'}) \}$. Hence the languages of the form
$\Phi(Forb(\mathcal{F'}))$ are closed under union. In the proof of
Theorem~\ref{NP} we may restrict ourselves to connected lifts when
proving that the constructed $\Phi(Forb(\mathcal{F'}))$ is the
desired language.
\end{remark}

Let us now formulate and prove the two analogous theorems for the
class monotone, monadic SNP and for the class monadic SNP without
inequality (which correspond to $(2)$ and $(3)$ of Theorem
\ref{123}). Here we use the categories $Rel^{cov}_{inj}(\Delta,
\Delta')$ and $Rel^{cov}_{full}(\Delta, \Delta')$.

\begin{theorem}~\label{inj}
For every language $L \in NP$ there exist relational types $\Delta$
and $\Delta'$, where $\Delta'$ contains only unary relational
symbols and a finite set $\mathcal{F'} \subset
Rel^{cov}_{inj}(\Delta, \Delta')$ such that $L$ is computationally
equivalent to the class $\Phi({\rm Forb}_{inj}(\mathcal{F'}))$.
\end{theorem}

\begin{proof}
We proceed analogously as in the proof of Theorem \ref{NP} for
formulas which are monotone monadic SNP. We stress the differences
only. First, using Lemma \ref{tech} again, we may suppose that $L$
is defined by a canonical primitive formula. This constitutes the
first step as now we do not have problem with the arity of the proof
relations since these are all monadic.

\noindent {\bf Step 2.}

\noindent
We want to enforce for $(\alpha_i \wedge \beta_i \wedge \varepsilon_i)$ and
distinct variables $x,y$ appearing in it that $x \neq y$ is an atom of
$\varepsilon_i$.
If this atom is not in $\beta_i$ then we exchange
$\neg (\alpha_i \wedge \beta_i \wedge \varepsilon_i)$ by the following
conjunction: $\neg (\alpha_{i_1} \wedge \beta_{i_1} \wedge \varepsilon_{i_1})
\bigwedge \neg (\alpha_{i_2} \wedge \beta_{i_2} \wedge \varepsilon_{i_2})$,
where $\neg (\alpha_{i_1} \wedge \beta_{i_1} \wedge \varepsilon_{i_1})$ is
$\neg (\alpha_i \wedge \beta_i \wedge \varepsilon_i)$ except that we replace
all occurence of $y$ by $x$ in it, $\alpha_{i_2}=\alpha_i,
\beta_{i_2}= \beta_i$ and $\varepsilon_{i_2}=\varepsilon_i \wedge (x\neq y)$.
This new formula is equivalent to the original one. In finitely many steps we
manage to enforce that all the required atoms of the form $x \neq y$ are there
in the appropriate $\varepsilon_i$.

\noindent

We now define $\Delta'$ in the same way as in Theorem~\ref{NP}
(thus $\Delta'$ is a monadic type). The set of forbidden lifts
$\mathcal{F'}$ is also defined analogously  as in Theorem~\ref{NP}
with the only one difference which relates to the construction of
formula $\varphi_{\bf F'}$ which will have now more clauses: the
formula $\varphi_{\bf F'}$ will have all the atom clauses as in
Theorem~\ref{NP} (i.e. $\varphi_{\bf F'}(x_1, \dots ,x_{|{\bf
F'}|})$ will contain as atoms all those tuples which  express the
fact that  a tuple $\ovl a$ is in the homomorphic image of $\bf
F'$) and additionally we will have atoms $x \neq y$ for every pair
of different variables. After this change we see easily that the
rest of the proof does not depend on which category we work.

\end{proof}

\begin{remark}~\label{pinj}
If we do not enforce the condition that the atom $x \neq y$
appears in every clause containing the variables $x$ and $y$ (Step
1. of the proof) before constructing $\mathcal{F'}$ then we get
some weaker characterization. Namely, the language $L$ will be
similar to the form of Theorem~\ref{inj} but we allow partially
injective mappings not only injective ones. For every ${\bf F'}
\in \mathcal{F'}$ and pair $x,y \in \bf F'$ we may have the plus
condition that they can not collapse by a homomorphism. The class
defined by such partially injective forbidden lifts still equals
to the class of languages of the form $\Phi({\rm
Forb}_{inj}(\mathcal{F'}))$: we can do Step 2. in this
combinatorial setting, too. Here the transformation means that for
any ${\bf F'} \in \mathcal{F'}$ and pair $x,y \in \bf F'$ which
may collapse, we exchange $\bf F'$ by two new forbidden
structures. One of the structures is $\bf F'$ with conditions on
the same pairs plus we require that $x$ and $y$ may not collapse.
The other is a factor of $\bf F'$ where we identify $x$ and $y$,
and we have the condition on a pair  of elements of new structures
not to collapse iff we have it on a pair in their preimages in
$\bf F'$. The iteration of this transformation expresses a
language defined in partially injective setting in the fully
injective terminology of Theorem~\ref{inj} (with the same $\Delta$
and $\Delta'$).
\end{remark}

\noindent
Let us now transform the third syntactic class of Theorem \ref{123}.

\begin{theorem}~\label{full}
For every language $L \in NP$ there exist relational types
$\Delta$ and $\Delta'$, where $\Delta'$ contains only unary
relational symbols and a finite set $\mathcal{F'} \subset
Rel^{cov}_{full}(\Delta, \Delta')$  such that $L$ is
computationally equivalent to the class $\Phi({\rm
Forb}_{full}(\mathcal{F'}))$.
\end{theorem}

\begin{proof}
The proof of Theorem \ref{full} is again a modification of  the
above proof of Theorem~\ref{NP} (and of Theorem~\ref{inj}) for
formulas in monadic SNP without inequality. The construction of
$\mathcal{F'}$ is even easier: Again, in Step 1., it suffices to
assume that $\varphi_L$ is canonical primitive. We only need to be
careful with construction of the formula $\varphi_{\bf F'}(x_1,
\dots ,x_{|{\bf F'}|})$ expressing the fact that the set $\{ x_1,
\dots ,x_{|\bf F'|} \}$ is the homomorphic image of $\bf F'$ (recall
that all homomorphism are now considered in $Rel_{full}(\Delta,
\Delta'))$. The formula will contain again more atoms. For every
tuple $\ovl{a}$ in the input relation $R$ we will have an atom
expressing that the image of the tuple is in relation $R$ like in
the proof of Theorem~\ref{NP}. But additionally we will have the
negation of such an atom for every tuple not contained by an input
relation. The rest of the proof is again the same.
\end{proof}

Similarly as above (Remark~\ref{pinj} to Theorem~\ref{inj})  we
have the possibility to state a weaker theorem in the notion of
partially full mappings. Consider a relational symbol $R$ of arity
$q$. We may have two conditions on an $q$-tuple in a structure
$\bf A$, either it is in $R$ or not. In the category $Rel_{full}$
this gives some restrictions on the homomorphisms of $\bf A$ in
both cases. We may generalize the class of objects such that for
every relation $R$ and $q$-tuple we have three possibilities (from
the viewpoint of mappings to a structure): either the tuple should
be mapped to a tuple in $R$, or to a tuple not in $R$ or we have
no restriction. We may define a class of languages in $Rel(\Delta)$
using this enlarged set of forbidden lifts. However this new class
of languages is still equal to those of the form $\Phi({\rm
Forb}_{full}(\mathcal{F'}))$. This may be seen as follows: a
forbidden lift in this new setting may be replaced by a set of
forbidden lift in $Rel^{cov}_{full}(\Delta, \Delta')$ as  for the
set of tuple-relation pairs without any condition we take all
possibilities of relation and non-relation conditions. This new
set of forbidden lifts defines the same language.

\section{Lifts and Shadows of Dualities}

Some of the transformations presented in Section 4 lead to deeper
results - the lifts and shadows give rise to a life on their own.
We prove here two results which will prove to be useful in the next section.

It follows from the Section 3 that shadows of classes ${\rm
Forb}(\mathcal F')$ (in three categories ${Rel^{cov}(\Delta,
\Delta')}$, ${Rel^{cov}_{inj}(\Delta, \Delta')}$ and
${Rel^{cov}_{full}(\Delta, \Delta')}$) include all NP-complete
languages. What about finite dualities? A delicate interplay of
lifting and shadows for dualities is expressed by the following two
statements which deal (for brevity) with classes $Rel^{cov}(\Delta,
\Delta')$ only. Despite its formal complexity Theorem \ref{eqdual1}
is an easy statement.

\smallskip

\begin{theorem}~\label{eqdual1}
Let ${\mathcal F'}$ be a finite set of structures in
${Rel^{cov}(\Delta, \Delta')}$. Suppose that there exist a finite
set of structures ${\mathcal D'}$ such that $({\mathcal
F'},{\mathcal D'})$ is a finite duality in ${Rel^{cov}(\Delta,
\Delta')}$. Then the following sets coincide: the shadow
$\Phi({\rm Forb}(\mathcal F')) = \{\Phi({\bf A'}):{\bf A'}\in {\rm
Forb}(\mathcal F')\}$ and $CSP(\Phi(\mathcal D')) = \bigvee_{{\bf
D'} \in \mathcal D'} CSP(\Phi({\bf D'}))$. Explicitly: for every
${\bf A } \in Rel(\Delta)$ there exists ${\bf A'} \in Rel(\Delta,
\Delta'), \Phi({\bf A'}) = \bf A$ and ${\bf F}'
\not\longrightarrow \bf A'$ for every ${\bf F}' \in \mathcal{F'}$
iff ${\bf A} \longrightarrow \Phi({\bf D'})$ for some ${\bf D}'
\in \mathcal{D'}$ .
\end{theorem}

Note that (in the statement of Theorem \ref{eqdual1}
we do not claim that the pair  $(\Phi({\mathcal
F'}),\Phi(\mathcal{D'}))$ is a duality in the class $Rel(\Delta)$.
This of course does not hold (as shown by our example of
$3$-colorability in the introduction). But the images of all
structures defined by all obstacles of $CSP({\mathcal D'})$ are
forming all obstacles of $CSP(\Phi({\mathcal D'}))$. We call this
{\it shadow duality}.

It is important that Theorem \ref{eqdual1} may be sometimes
reversed: shadow dualities may be sometimes ``lifted''. This is
non-trivial and in fact Theorem \ref{eqdual2} may be seen as the
core of this paper.

\noindent
\begin{theorem}~\label{eqdual2}
Let ${\mathcal F'}$ be a finite set of structures in
${Rel^{cov}(\Delta, \Delta')}$, consider ${\rm Forb}(\mathcal F')$
and suppose that $\Phi({\rm Forb}(\mathcal F')) = CSP({\mathcal
D})$ (in $Rel(\Delta)$) for a finite set ${\mathcal D}$ of objects
of $Rel(\Delta)$. (In other words let the pair $(\mathcal
F',\mathcal D)$ form a shadow duality.) Assume also that
$CSP({\mathcal D}) \neq Rel(\Delta)$ and that $\Delta'$ contains
only unary relations. Then there exists a finite set ${\mathcal
D'}$ in ${Rel^{cov}(\Delta, \Delta')}$ such that ${\rm
Forb}(\mathcal F') = CSP({\mathcal D'})$.
\end{theorem}

\noindent
Before proving Theorems~\ref{eqdual1} and
\ref{eqdual2} we formulate first a lemma which we shall use repeatedly:

\begin{lemma}(lifting)
\label{lifting} Let ${\bf A, B} \in Rel(\Delta)$, homomorphism $f:
{\bf A} \longrightarrow \bf B$ and  $\Phi(\bf B') = \bf B$ be
given. Then there exists ${\bf A'} \in Rel^{cov}(\Delta,
\Delta')$, $\Phi({\bf A'})={\bf A}$ such that the mapping $f$ is a
homomorphism ${\bf A'} \longrightarrow \bf B'$ in
$Rel^{cov}(\Delta, \Delta')$.
\end{lemma}

\begin{figure}[!h]
\begin{center}
\includegraphics[width=4.5cm]{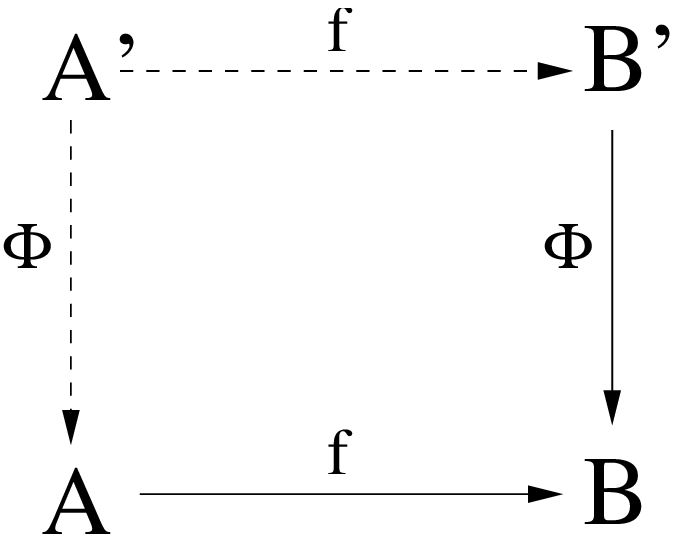}
\end{center}
\end{figure}

\vspace{10mm}

\begin{proof}
Assume that ${\bf A, B} \in Rel(\Delta)$, $\Phi(\bf{B'}) = \bf B$ and
$f: {\bf A} \longrightarrow \bf B$ are as in the statement. For
each $R \in \Delta'$ put $R({\bf A'}) = f^{-1}(R(\bf{B'}))$. It is
easy to see that ${\bf A'} \in Rel^{cov}(\Delta,\Delta')$
\end{proof}

\begin{proof}(of Theorem \ref{eqdual1})
Suppose that ${\bf A} \in CSP(\Phi(\mathcal D'))$, say ${\bf A}
\in CSP(\Phi({\bf D}'))$. Now for a homomorphism $f: {\bf
A}\longrightarrow \Phi({\bf D}') $ there is at least one lift
${\bf A}'$ of $\bf A$ such that the mapping $f$ is a homomorphism
${\bf A}' \rightarrow {\bf D}'$ (here we use Lifting Lemma
\ref{lifting}). By the duality $({\mathcal F'},{\mathcal D'})$ (in
$Rel^{cov}(\Delta, \Delta')$) ${\bf F}' \nrightarrow {\bf A'}$ for
any ${\bf F}' \in \mathcal D'$ and thus in turn ${\bf A} \in
\Phi({\rm Forb}(\mathcal F'))$.

Conversely, let us assume  that ${\bf A'} \in {\rm Forb}(\mathcal
F')$ satisfies $\Phi(\bf A') = \bf A$. But then ${\bf A'} \in
CSP(\mathcal D')$ and thus by the functoral property of $\Phi$ we
have ${\bf A} = \Phi({\bf A'}) \in CSP(\Phi(\mathcal D'))$.

\end{proof}

\begin{proof}(of Theorem \ref{eqdual2})
Assume $\Phi({\rm Forb}(\mathcal F')) = CSP({\mathcal D})$. Our goal
is to find $\mathcal D'$ such that ${\rm Forb}(\mathcal F') =
CSP({\mathcal D'})$. This will follow as a (non-trivial) combination
of Theorem \ref{FNT} and Lemma \ref{sparse}. By Theorem \ref{FNT} we
know that if $\mathcal F'$ is a set of (relational) forests then the
set $\mathcal F'$ has a dual set $\mathcal D'$ (in the class
$Rel^{cov}(\Delta,\Delta')$). So assume to the contrary that one of
the structures, say ${\bf F}_0'$, fails to be a forest (i.e. we
assume that one of the components of ${\bf F}_0'$ has a cycle). We
proceed by a refined induction (which will allow us to use more
properties of ${\bf F}_0'$). Let us introduce carefully the setting
of the induction.

We assume shadow duality $\Phi({\rm Forb}(\mathcal F')) =
CSP({\mathcal D})$. Let ${\mathcal D}$ be fixed throughout the
proof. Clearly many sets $\mathcal F'$ will do the job and we
select the set $\mathcal F'$ such that $\mathcal F'$ consists of
cores of all homomorphic images (explicitly: we close $\mathcal
F'$ on homomorphic images and then take the set of  cores of all
these structures). Among all such sets  $\mathcal F'$ we take a
set of  minimal cardinality. It will be again denoted by $\mathcal
F'$. We proceed by induction on the size $|\mathcal F'|$ of
$\mathcal F'$.

The set ${\rm Forb}(\mathcal F')$ is clearly determined by the
minimal elements of $\mathcal F'$ (minimal in the homomorphism
order). Thus let us assume that one of these minimal elements, say
${\bf F}_0'$, is not a forest. By the minimality of $\mathcal F'$
we see that we have a proper inclusion $\Phi({\rm Forb}(\mathcal
F' \setminus \{{\bf F}_0'\})) \supset CSP({\mathcal D})$. Thus
there exists a structure $\bf S$ in the difference. But this in
turn means that there has to be a lift $\bf S'$ of $\bf S$ such
that ${\bf F}_0' \longrightarrow \bf S'$ and ${\bf S}
\not\rightarrow {\bf D}$ for every ${\bf D} \in \mathcal D$. In
fact not only that: as ${\bf F}_0'$ is a core, as ${\rm
Forb}(\mathcal F')$ is homomorphism closed  and as $\mathcal F'$
has minimal size we conclude that there exist $\bf S$ and $\bf S'$
such that any homomorphism ${\bf F}_0' \longrightarrow {\bf S}'$
is a monomorphism (i.e. one-to-one, otherwise we could replace
${\bf F}_0'$ by a set of all its homomorphic images - ${\bf F}_0'$
would not be needed).

Now we apply (the second non-trivial ingredient) Lemma \ref{sparse}
to the structure $\bf S$ and an $\ell > |X({\bf F}_0')|$: we find a
structure ${\bf S}_0$ with the following properties: ${\bf S}_0
\longrightarrow \bf S$,  ${\bf S}_0 \longrightarrow {\bf D}$ if and
only if  ${\bf S} \longrightarrow {\bf D}$ for every ${\bf D} \in
{\mathcal D}$ and ${\bf S}_0$ contains no cycles of length $\leq
\ell$. It follows that  ${\bf S}_0 \not\in CSP({\mathcal D})$. Next
we apply Lemma \ref{lifting} to obtain a structure ${\bf S}_0'$ with
${\bf S}_0' \longrightarrow \bf S'$. Now we use that all relations
in $\Delta'$ are unary and we see that ${\bf S}_0'$ does not contain
cycles of length $\leq \ell$. Now for any ${\bf F}' \in \mathcal
F'$, ${\bf F}' \neq {\bf F}_0'$ we have ${\bf F}' \nrightarrow {\bf
S}_0'$ as ${\bf S}_0' \rightarrow \bf S'$ and ${\bf F}' \nrightarrow
{\bf S}'$. As the only homomorphism ${\bf F}_0' \longrightarrow {\bf
S}'$ is a monomorphism the only (hypothetical) homomorphism ${\bf
F}_0' \longrightarrow \bf S'$ is also monomorphism. But this is a
contradiction as ${\bf F}_0'$ contains a cycle while ${\bf S}_0'$
has no cycles of length $\leq \ell$. This completes the proof.
\end{proof}

\section{MMSNP and forbidden patterns}

Madelaine \cite{M} introduced the class FP. Every language of the
class FP is defined by {\it forbidden patterns} which are defined as
follows. Consider the finite relational type $\Delta$, the finite
set $T$ and the set of pairs $({\bf F}_1, \varphi_1), \dots ,({\bf
F}_n, \varphi_n)$, where each ${\bf F}_i \in Rel(\Delta)$ and
$\varphi_i: {\bf F}_i \rightarrow T$ is a mapping $(i = 1, \ldots,
n)$. The language $L$ belongs to the class FP if there are patterns
$({\bf F}_1, \varphi_1), \dots ,({\bf F}_n, \varphi_n)$ such that
$L$ is the class of all structures ${\bf A} \in Rel(\Delta)$  for
which there exists a mapping $\varphi: {A} \rightarrow T$ such that
for all $i = 1,\ldots, n$  no homomorphism $\alpha: {\bf F}_i
\rightarrow {\bf A}$ satisfies $\varphi \circ \alpha \neq
\varphi_i$. Formally: $L=\{{\bf A} \in Rel(\Delta): \exists \varphi:
{\bf A} \rightarrow T \text{ such that } \forall i, \alpha: {\bf
F}_i \rightarrow {\bf A} \text{ homomorphism }  \varphi \circ \alpha
\neq \varphi_i \}$.

This is a special case of our approach and the class FP may
be equivalently defined as follows (using lifts and shadows):  we say that
the set $L \subseteq Rel(\Delta)$ is an {\em FP-language} if there exist a finite type $\Delta'$
of monadic (unary) relational symbols and a language $L' \in Rel^{cov}(\Delta, \Delta')$ such that
$L = \Phi(Forb(\mathcal{F'}))$ for a finite set $\mathcal{F'} \subseteq
Rel^{cov}(\Delta, \Delta')$. (Thus $\Delta'$ is a partition on every ${\bf F}' \in \mathcal{F'}$.)
The equivalence  is clear: we
consider the signature (relational type) $\Delta'$ that contains the
unary symbol  $u_t$ for every element $t \in T$.
To every pattern $({\bf F}_i,\varphi_i)$ we correspond the
relational structure ${\bf F}'_i \in Rel(\Delta, \Delta')$ with the
shadow
${\bf F}_i$ such that the element $x \in {\bf F}_i$ is in the relation
$u_{\varphi_i(x)}$. The converse is also evident: every
FP-language can be defined by forbidden patterns.

In other words the class FP is the class of languages defined by forbidden
monadic lifts of the class $Rel(\Delta)$.

It has been proved  in \cite{M} that the classes FP and MMSNP are
equal. The containment FP $\supseteq$ MMSNP follows from the proof
of Theorem \ref{NP}: every MMSNP problem (as any NP problem) can be
considered as the shadow of a language ${\rm Forb}({\mathcal F'})$
in a lifted category $Rel^{cov}(\Delta, \Delta')$. In the case of
the class MMSNP these lifted relations (in $\Delta'$) are all unary.
And for unary relations we use the preceding remark which claims
that monadic lifts and forbidden patterns are equivalent. In order
to prove the converse one needs to show that every language defined
by forbidden monadic lifts is in MMSNP. This part of proof is
straightforward.

Madelaine and Stewart \cite{MStuart} gave a long process to decide
whether an FP language is a finite union of CSP languages. We use
Theorems~\ref{eqdual1}, \ref{eqdual2} and the description of
dualities for relational structures \cite{FNT} to give a short
characterization of a more general class of  languages.

\begin{theorem}
\label{main}
Consider the language $L$ determined by forbidden monadic lifts. Explicitly,
$L = \Phi(Forb(\mathcal{F'}))$ for a set ${\mathcal F'} \subset
Rel(\Delta, \Delta')$ (with $\Delta'$ monadic).
If no ${\bf F'} \in Forb(\mathcal{F'})$ contains a cycle then there is a
set of finite structures $\mathcal{D} \subseteq Rel(\Delta)$ such that
$L = CSP(\mathcal{D})$. If one of the lifts ${\bf F'}$ in a minimal
subfamily of $\mathcal{F'}$ contains a cycle in its core then the
language $L$ is not a finite union of CSP languages.
\end{theorem}

\begin{proof}
If no ${\bf F}' \in Forb(\mathcal{F'})$ contains a cycle then the
set  $\mathcal F'$ has a dual $\mathcal D'$  in
${Rel^{cov}(\Delta, \Delta')}$ by \cite{FNT}, and the shadow of
this set $\mathcal D'$  gives the dual set $\mathcal{D}$ of the
set $\Phi(Forb(\mathcal{F'}))$ (by Theorem~\ref{eqdual1}). On the
other side  if one ${\bf F}' \in Forb(\mathcal{F'})$ contains a
cycle in its core and if $\mathcal{F'}$ is minimal (i.e. $\bf F'$
is needed) then $Forb(\mathcal{F'})$ does not have a dual in
$Rel^{cov}(\Delta, \Delta')$. The shadow of the language
$Forb(\mathcal{F'})$ is the language $L$ and consequently this
fails to be a finite union of CSP languages by
Theorem~\ref{eqdual2} (as every monadic shadow duality can be
lifted).
\end{proof}

\section{Understanding Feder - Vardi}

Now we prove one of the principal results of \cite{FV}  by tools
which we developed in previous sections. Feder and Vardi have proved
that the classes MMSNP and CSP are random equivalent, this was later
derandomised. Here we discuss the deterministic part of the
Feder-Vardi proof. It seems that our setting streamlines some of the
earlier arguments. Our proof is not essentially different from the
original one, yet the use of dualities makes the construction of the
proof natural and easier.

A structure ${\bf A}$ is {\it biconnected} if
every point deleted substructure is connected (in other words
for every three distinct elements $x,y$ and $z$ there is a path
connecting $x$ and $y$ that avoids $z$). Note that a biconnected
structure with more than one relational tuple contains a cycle.
Inclusion maximal biconnected substructures are called {\it
biconnected components} (in graph theory they are called
blocks). For the set of relational structures $\mathcal{D}$ we
denote by $CSP_{girth>k}(\mathcal{D})$ the language of structures in
$CSP(\mathcal{D})$ with girth larger than $k$. We will prove the
following:

\begin{theorem}\cite{FV}
\label{FV} For every MMSNP language $L$ there is a finite set of
relational structures $\mathcal{D}$ (of possibly different type) and
a positive integer $k$ such that the following hold.

\begin{enumerate}

\item $L$ can be polynomially reduced to $CSP(\mathcal{D})$.

\item The language $CSP_{girth>k}(\mathcal{D})$ can be polynomially reduced to $L$.
\end{enumerate}
\end{theorem}

\begin{proof}
We assume that $L = \Phi(Forb(\mathcal{F'}))$ for a set ${\mathcal F'} \subset
Rel(\Delta, \Delta')$ (with $\Delta'$ monadic).
We construct the set $\mathcal{D}$. First we determine the
type of the relational structures in $\mathcal{D}$.

Let $\mathcal{B}
\subset Rel(\Delta,\Delta')$ be the set of the biconnected
components of the structures in $\mathcal{F'}$. For every
isomorphism class in $\Phi(\mathcal{B})$ we choose one structure ${\bf
B}$ in this class. We may assume that the
base set of the representative ${\bf B}$ is $\{1, \dots ,|{\bf B}|
\}$. For each ${\bf B} \in \Phi(\mathcal{B})$
we introduce the relational symbol $R_{\bf B}$ of arity  $|{\bf B}|$
(the size of the structure). We denote
by $\beta$ the type that consists of these relational symbols
$R_{{\bf B}}$. This will be the type of the structures in
$\mathcal{D}$.

Next, we define the following functors $\Psi$ and $\Theta$. The functor
$\Psi: Rel(\Delta) \rightarrow Rel(\beta)$ assigns to a structure
${\bf A}$ a structure $\Psi({\bf A})$. $\bf A$ and $\Psi({\bf A})$ have the same base set
and its relatons are defined as follows:

$R_{\bf B}(\Psi({\bf A}))= \{ f | f: {\bf B} \rightarrow {\bf A}
\}$ ($f$ is a homomorphism in $Rel(\Delta)$).

I.e. a tuple of elements is in $R_{\bf B}$ relation if it is the
homomorphic image of ${\bf B}$. The functor $\Theta$ maps to a
structure ${\bf A} \in Rel(\beta)$ the following structure again on the
same base set in $Rel(\Delta)$:

$\Theta({\bf A})= \cup \{ f({\bf B}): {\bf B} \in \Phi(\mathcal{B}),
f \in R_{\bf B}({\bf A}) \}$ (here $f({\bf B})$ is the homomorphic image of the structure
$\bf B$).

The mappings $\Psi$ and $\Theta$ are both functoral. Consider the
induced functors $\Psi': Rel(\beta,\Delta') \rightarrow
Rel(\Delta,\Delta')$ and $\Theta': Rel(\Delta,\Delta') \rightarrow
Rel(\beta, \Delta')$. We will use the following properties of these
functors.

(i) $\Theta \circ \Psi = id_{Rel(\Delta)}$ and $\Theta' \circ \Psi'
= id_{Rel(\Delta,\Delta')}$

(ii) For every ${\bf B} \in Rel(\beta)$ (${\bf B'} \in
Rel(\beta,\Delta')$) and relational symbol $R \in \beta$ the
following inclusions hold:

 $R(\Psi \circ \Theta ({\bf B})) \supseteq R({\bf B})$
$\big( R(\Psi' \circ \Theta' ({\bf B'})) \supseteq R({\bf B'})
\big)$.

We continue the construction of $\mathcal{D}$. We define the
finite set of structures $\mathcal{G'} \subset Rel(\beta,
\Delta')$ as follows: We put ${\bf G'} \in \mathcal{G'}$ if

\begin{enumerate}
\item
${\bf G'}$ is a forest,

\item
there exists ${\bf F'} \in \mathcal{F'}$ and a homomorphism $\varphi:
{\bf F'} \rightarrow \ \Theta({\bf G'})$, such that
every element of ${\bf G'}$ is contained in $\varphi({\bf F'})$ or in a
relational tuple intersecting $\varphi({\bf F'})$.

\end{enumerate}

Observe the following straightforward consequences of the
construction of $\mathcal{G'}$.

(iii) If ${\bf A'} \notin Forb(\mathcal{F'})$ then $\Psi({\bf A'})
\notin Forb(\mathcal{G'})$ holds for every ${\bf A'} \in
Rel(\Delta,\Delta')$, since ${\bf A'} \in \mathcal{F'}
\Longrightarrow \Psi({\bf A'}) \notin Forb(\mathcal{G'})$.

(iv) If ${\bf B'} \notin Forb(\mathcal{G'})$ then $\Theta({\bf B'})
\notin Forb(\mathcal{F'})$ holds for every ${\bf B'} \in
Rel(\beta,\Delta')$, since ${\bf B'} \in \mathcal{G'}
\Longrightarrow \Theta({\bf B'}) \notin Forb(\mathcal{F'})$.

The set $\mathcal{G'}$ consists of finitely many relational forests.
Hence we know by Theorem~\ref{FNT} that
$\Phi(Forb(\mathcal{G'}))=CSP(\mathcal{D})$ for some finite
$\mathcal{D} \subset Rel(\beta)$. We will prove that the conditions
of the theorem hold for this choice of $\mathcal{D}$. All the
reductions will be functoral.

First we prove that $L$ can be polynomially reduced to
$\Phi(Forb(\mathcal{G'}))$. We succeed to show that for a structure
${\bf A} \in Rel(\Delta)$ the equivalence ${\bf A} \in L \iff
\Psi({\bf A}) \in \Phi(Forb(\mathcal{G'}))$ holds. This is
implied by the equivalence in the lifted category as
the same $\Delta'$ relations prove the membership in both
languages: If ${\bf A'} \in Forb(\mathcal{F'})$ then $\Psi'({\bf
A'}) \in Forb(\mathcal{G'})$ by (i) and (iv). On the other hand
(iii) implies that if ${\bf A'} \notin Forb(\mathcal{F'})$ then
$\Psi'({\bf A'}) \notin Forb(\mathcal{G'})$.

Let $k$ denote the size of the largest structure in $\mathcal{F'}$.
We prove that $CSP_{girth>k}(\mathcal{D})$ can be polynomially
reduced to $L$. In fact we will prove that for every ${\bf B} \in
Rel(\beta)$ with girth $>k$ the equivalence ${\bf B} \in
\Phi(Forb(\mathcal{G'})) \iff \Theta({\bf B}) \in L$ holds. Again we
prove the equivalence in the lifted categories. If
$\Theta'({\bf B'}) \in Forb(\mathcal{G'})$ then $\Psi'(\Theta'({\bf
B'})) \in Forb(\mathcal{F'})$, as we have seen in the reduction of
$L$ to $K$. The structure ${\bf B'}$ contains less relations than
$\Psi'(\Theta'({\bf B'}))$ by (ii), hence ${\bf B'} \in
Forb(\mathcal{G'})$. If $\Theta'({\bf B'}) \notin
Forb(\mathcal{G'})$ then there exists a structure ${\bf F'} \in
\mathcal{F'}$ such that $\varphi:{\bf F'} \rightarrow {\bf B'}$. By
the girth condition on ${\bf B'}$ we know that the union of
$\varphi({\bf F'})$ and the relational tuples intersecting
$\varphi({\bf F'})$ is a forest. Hence there is a structure ${\bf
G'} \in \mathcal{G'}$ isomorphic to this substructure of ${\bf B'}$.
Now ${\bf G'} \rightarrow {\bf A'}$, hence ${\bf A'} \notin
Forb(\mathcal{G'})$.
\end{proof}

The remaining part is the reduction of CSP with large girth to CSP. Feder
and Vardi proved a randomized reduction, this was later derandomized.

\begin{lemma}[\cite{CSPvsMMSNP}]
For every finite set of relational structures $\mathcal{D}$
and integer $k>0$ the language $CSP(\mathcal{D})$ can be
polynomially reduced to $CSP_{girth>k}(\mathcal{D})$.
\end{lemma}

The essence of this reduction is the Sparse Incomparability Lemma \ref{sparse}.
This polynomial reduction was proved with expanders in the case of digraphs
\cite{MN}. The reduction in the case of general relational structures needed a
generalization of expanders called expander (relational) structures. The
notion of expander relational structures was introduced in \cite{CSPvsMMSNP}
\cite{K}, and also a polynomial time construction of such structures with
large girth is given there.

\end{document}